# An endpoint estimate for the Kunze-Stein phenomenon and related maximal operators

By Alexandru D. Ionescu*

## Abstract

One of the purposes of this paper is to prove that if $G$ is a noncompact connected semisimple Lie group of real rank one with finite center, then

$$L^{2,1}(G) * L^{2,1}(G) \subseteq L^{2,\infty}(G).$$

Let $K$ be a maximal compact subgroup of $G$ and $X = G/K$ a symmetric space of real rank one. We will also prove that the noncentered maximal operator

$$\mathcal{M}_2 f(z) = \sup_{z \in B} \frac{1}{|B|} \int_B |f(z')| \, dz'$$

is bounded from $L^{2,1}(X)$ to $L^{2,\infty}(X)$ and from $L^p(X)$ to $L^p(X)$ in the sharp range of exponents $p \in (2, \infty]$. The supremum in the definition of $\mathcal{M}_2 f(z)$ is taken over all balls containing the point $z$.

## 1. Introduction

A central result in the theory of convolution operators on semisimple Lie groups is the Kunze-Stein phenomenon which, in its classical form, states that if $G$ is a connected semisimple Lie group with finite center and $p \in [1, 2)$, then

(1.1) $$L^2(G) * L^p(G) \subseteq L^2(G).$$

The usual convention, which will be used throughout this paper, is that if $\mathcal{U}$, $\mathcal{V}$, and $\mathcal{W}$ are Banach spaces of functions on $G$ then the notation $\mathcal{U} * \mathcal{V} \subseteq \mathcal{W}$ indicates both the set inclusion and the associated norm inequality. The inclusion (1.1) was established by Kunze and Stein [10] in the case when the group $G$ is $\mathrm{SL}(2, \mathbb{R})$ (and, later on, for a number of other particular groups) and by Cowling [3] in the general case stated above. For a more complete account of the development of ideas leading to (1.1) we refer the reader to [3] and [4].

*The author was supported by an Alfred P. Sloan graduate fellowship.

More recently, Cowling, Meda and Setti noticed that if the group $G$ has real rank one then the inclusion (1.1) can be strengthened. Following earlier work of Lohoué and Rychener [9], the key ingredient in their approach is the use of Lorentz spaces $L^{p,q}(G)$; they prove in [4] that if $G$ is a connected semisimple Lie group of real rank one with finite center, $p \in (1,2)$ and $(u,v,w) \in [1,\infty]^3$ has the property that $1 + 1/w \leq 1/u + 1/v$, then

$$(1.2) \qquad L^{p,u}(G) * L^{p,v}(G) \subseteq L^{p,w}(G).$$

In particular, $L^{p,1}$ convolves $L^p$ into $L^p$ for any $p \in [1,2)$. Our first theorem is an endpoint estimate for (1.2) showing what happens when $p = 2$.

THEOREM A. *If $G$ is a noncompact connected semisimple Lie group of real rank one with finite center then*

$$(1.3) \qquad L^{2,1}(G) * L^{2,1}(G) \subseteq L^{2,\infty}(G).$$

Notice that (1.2) follows from Theorem A and a bilinear interpolation theorem ([4, Theorem 1.2]). Unlike the classical proofs of the Kunze-Stein phenomenon, our proof of Theorem A will be based on real-variable techniques only: the inclusion (1.3) is equivalent to an inequality involving a triple integral on $G$ and we use certain nonincreasing rearrangements to control this triple integral. Easy examples, involving only $K$-bi-invariant functions, show that the inclusion (1.3) is sharp in the sense that neither of the $L^{2,1}$ spaces nor the $L^{2,\infty}$ space can be replaced with some $L^{2,u}$ space for any $u \in (1,\infty)$.

Let $K$ be a maximal compact subgroup of the group $G$ and $X = G/K$ the associated symmetric space. Assume from now on that the group $G$ satisfies the hypothesis stated in Theorem A and let $d$ be the distance function on $X \times X$ induced by the Killing form on the Lie algebra of the group $G$. Let $B(x,r)$ denote the ball in $X$ centered at the point $x$ of radius $r$ (with respect to the distance function $d$) and let $|A|$ denote the measure of the set $A \subset X$. For any locally integrable function $f$ on $X$, let

$$(1.4) \qquad \mathcal{M}_2 f(z) = \sup_{z \in B} \frac{1}{|B|} \int_B |f(z')|\, dz',$$

where the supremum in the definition of $\mathcal{M}_2 f(z)$ is taken over all balls $B$ containing $z$. We will prove the following:

THEOREM B. *The operator $\mathcal{M}_2$ is bounded from $L^{2,1}(X)$ to $L^{2,\infty}(X)$ and from $L^p(X)$ to $L^p(X)$ in the sharp range of exponents $p \in (2,\infty]$.*

We recall that the more standard centered maximal operator

$$\mathcal{M}_1 f(z) = \sup_{r>0} \frac{1}{|B(z,r)|} \int_{B(z,r)} |f(z')|\, dz'$$

is bounded from $L^1(X)$ to $L^{1,\infty}(X)$ and from $L^p(X)$ to $L^p(X)$ for any $p > 1$, as shown in [5] and [12] (without the assumption that $G$ has real rank one). Notice however that, unlike in the case of Euclidean spaces, balls on symmetric spaces do not have the basic doubling property (i.e. $|B(z, 2r)|$ is not proportional to $|B(z, r)|$ if $r$ is large), thus the maximal operators $\mathcal{M}_1$ and $\mathcal{M}_2$ are not comparable. Easy examples (see [7, Section 4]) show that Theorem B is sharp in the sense that the maximal operator $\mathcal{M}_2$ is not bounded from $L^{2,u}(X)$ to $L^{2,v}(X)$ unless $u = 1$ and $v = \infty$.

This paper is organized as follows: in the next section we recall most of the notation related to semisimple Lie groups and symmetric spaces and prove a proposition that explains the role of the Lorentz space $L^{2,1}(G//K)$ – the subspace of $K$-bi-invariant functions in $L^{2,1}(G)$. In Section 3 we prove Theorem B. As a consequence of Theorem B we obtain in Section 4 a covering lemma on noncompact symmetric spaces of real rank one. In Section 5 we give a complete proof of Theorem A, which is divided into four steps. The main estimate in the proof of Theorem A uses the technique of nonincreasing rearrangements; we return to this technique in the last section and prove a general rearrangement inequality.

We conclude this section with some remarks on semisimple Lie groups of higher real rank. If the group $G$ has real rank different from 1, then (1.2) fails (the estimate in Lemma 6 and the discussion following Proposition 7 in [1] show that the appropriate spherical function $\Phi_p$ fails to belong to $L^{p',\infty}(G)$, where $p'$ is the conjugate exponent of $p$); therefore Theorem A fails to hold. On the other hand, the author has recently proved by a different method in [7] that the $L^p$ estimate in Theorem B holds on symmetric spaces of arbitrary real rank. In the general case it is not known however whether the maximal operator $\mathcal{M}_2$ is bounded from $L^{2,1}(X)$ to $L^{2,\infty}(X)$.

This work is part of the author's doctoral thesis at Princeton University under the guidance of Prof. Elias M. Stein. I would like to thank him for many clarifying discussions on the subject and for his time, interest and steady support. I would also like to thank Jean-Philippe Anker for several corrections on a preliminary version of this paper and the referee of the paper for a careful and detailed report.

## 2. Preliminaries

Let $G$ be a noncompact connected semisimple Lie group with finite center, and let $\mathfrak{g}$ be its Lie algebra. Most of our notation related to semisimple Lie groups and symmetric spaces is standard and can be found for example in [6]. Fix a Cartan involution $\theta$ of $\mathfrak{g}$ and let $\mathfrak{g} = \mathfrak{k} \oplus \mathfrak{p}$ be the associated Cartan decomposition. Let $\mathfrak{a}$ be a maximal abelian subspace of $\mathfrak{p}$; we will assume from now on that the group $G$ has real rank one, i.e., $\dim \mathfrak{a} = 1$. Let $\mathfrak{a}^*$ denote the real dual of $\mathfrak{a}$, let $\Sigma \subset \mathfrak{a}^*$ be the set of nonzero roots of the pair $(\mathfrak{g}, \mathfrak{a})$ and let $W$ be the Weyl group associated to $\Sigma$. It is well-known that $W = \{1, -1\}$ and $\Sigma$ is either of the form $\{-\alpha, \alpha\}$ or of the form $\{-2\alpha, -\alpha, \alpha, 2\alpha\}$. Let $m_1 = \dim \mathfrak{g}_{-\alpha}$, $m_2 = \dim \mathfrak{g}_{-2\alpha}$, $\rho = \frac{1}{2}(m_1 + 2m_2)\alpha$ and $\mathfrak{a}_+ = \{H \in \mathfrak{a} : \alpha(H) > 0\}$. Finally let $\overline{\mathfrak{n}} = \mathfrak{g}_{-\alpha} + \mathfrak{g}_{-2\alpha}$, $\overline{N} = \exp \overline{\mathfrak{n}}$, $K = \exp \mathfrak{k}$, $A = \exp \mathfrak{a}$ and $A_+ = \exp \mathfrak{a}_+$ and let $X = G/K$ be a symmetric space of real rank one.

The group $G$ has an Iwasawa decomposition $G = \overline{N}AK$ and a Cartan decomposition $G = K\overline{A_+}K$. Our proofs are based on relating these two decompositions, and for real rank one groups one has the explicit formula in [6, Ch.2, Theorem 6.1]. A similar idea was used by Strömberg [12] for groups of arbitrary real rank. Let $H_0 \in \mathfrak{a}$ be the unique element of $\mathfrak{a}$ for which $\alpha(H_0) = 1$ and let $a(s) = \exp(sH_0)$ for $s \in \mathbb{R}$ be a parametrization of the subgroup $A$. By [6, Ch.2, Theorem 6.1] one can identify the group $\overline{N}$ with $\mathbb{R}^{m_1} \times \mathbb{R}^{m_2}$ using a diffeomorphism $\overline{n} : \mathbb{R}^{m_1} \times \mathbb{R}^{m_2} \to \overline{N}$. This diffeomorphism has the property that if $t \geq 0$ then $\overline{n}(v, w)a(s) \in Ka(t)K$ if and only if

$$(2.1) \qquad (\cosh t)^2 = \left[\cosh s + e^s |v|^2\right]^2 + e^{2s}|w|^2.$$

In addition,

$$(2.2) \qquad a(s)\overline{n}(v, w)a(-s) = \overline{n}(e^{-s}v, e^{-2s}w).$$

Let $|\rho| = \rho(H_0) = \frac{1}{2}(m_1 + 2m_2)$ and let $dg$, $d\overline{n}$ and $dk$ denote Haar measures on $G$, $\overline{N}$ and $K$, the last one normalized such that $\int_K 1\, dk = 1$. Then the following integral formulae hold for any continuous function $f$ with compact support:

$$(2.3) \qquad \int_G f(g)\, dg = C_1 \int_K \int_{\mathbb{R}_+} \int_K f(k_1 a(t) k_2)(\sinh t)^{m_1}(\sinh 2t)^{m_2}\, dk_2\, dt\, dk_1,$$

and

$$(2.4) \qquad \int_G f(g)\, dg = C_2 \int_K \int_{\mathbb{R}} \int_{\overline{N}} f(\overline{n}a(s)k) e^{2|\rho|s}\, d\overline{n}\, ds\, dk$$
$$= C_2' \int_K \int_{\mathbb{R}} \int_{\mathbb{R}^{m_1} \times \mathbb{R}^{m_2}} f(\overline{n}(v, w)a(s)k) e^{2|\rho|s}\, dv\, dw\, ds\, dk.$$

The measures $dv$ and $dw$ are the usual Lebesgue measures on $\mathbb{R}^{m_1}$ and $\mathbb{R}^{m_2}$, and the constants $C_1$, $C_2$ and $C_2'$ depend on the normalizations of the various Haar measures. We will need a new integration formula, which is the subject of the following lemma.

LEMMA 1. *Suppose that $f : G \to \mathbb{C}$ is a $K$-bi-invariant (i.e., $f(k_1 g k_2) = f(g)$ for any $k_1, k_2 \in K$) continuous function with compact support and $F(t) = f(a(t))$ for any $t \in [0, \infty)$. Then for any $s \in \mathbb{R}$*

$$e^{|\rho|s} \int_{\overline{N}} f(\overline{n}a(s))\, d\overline{n} = \int_{|s|}^{\infty} F(t)\psi(t, s)\, dt,$$

*where the kernel $\psi : \mathbb{R}_+ \times \mathbb{R} \to \mathbb{R}_+$ has the property that $\psi(t,s) = 0$ if $t < |s|$ and*

$$(2.5) \qquad \psi(t,s) \approx \sinh t (\cosh t)^{m_2/2} (\cosh t - \cosh s)^{(m_1+m_2-2)/2}$$

*if $t \geq |s|$.*

As usual, the notation $U \approx V$ means that there is a constant $C \geq 1$ depending only on the group $G$ such that $C^{-1} U \leq V \leq CU$. This lemma is essentially proved in [8, Section 5]. For later reference we reproduce its proof.

*Proof of Lemma* 1. For any $t \geq |s|$, let

$$(2.6) \quad T_{t,s} = \{(v,w) \in \mathbb{R}^{m_1} \times \mathbb{R}^{m_2} : (\cosh t)^2 = \left[\cosh s + e^s |v|^2\right]^2 + e^{2s}|w|^2\}$$

be the set of points $P = P(v,w) \in \mathbb{R}^{m_1} \times \mathbb{R}^{m_2}$ with the property that $\overline{n}(P)a(s) \in Ka(t)K$ (these surfaces will play a key role in the proof of Theorem A). Let $d\omega_{t,s}$ be the induced measure on $T_{t,s}$ such that

$$\int_{\mathbb{R}^{m_1} \times \mathbb{R}^{m_2}} \phi(v,w) \, dv \, dw = \int_{t \geq |s|} \left[\int_{T_{t,s}} \phi(P) d\omega_{t,s}(P)\right] dt$$

for any continuous compactly supported function $\phi$. Then, since the function $f$ is $K$-bi-invariant,

$$e^{|\rho|s} \int_N f(\overline{n}a(s)) \, d\overline{n} = C e^{|\rho|s} \int_{\mathbb{R}^{m_1} \times \mathbb{R}^{m_2}} f(\overline{n}(v,w)a(s)) \, dv \, dw$$

$$= C e^{|\rho|s} \int_{t \geq |s|} F(t) \left[\int_{T_{t,s}} 1 \, d\omega_{t,s}\right] dt.$$

Let $\psi(t,s) = e^{|\rho|s} \int_{T_{t,s}} 1 \, d\omega_{t,s}$ and assume that $m_2 \geq 1$. We make the change of variables $v = [e^{-s}(u \cosh t - \cosh s)]^{1/2} \omega_1$ and $w = e^{-s} \cosh t (1-u^2)^{1/2} \omega_2$, where $\omega_1 \in S^{m_1-1}$ (the $m_1 - 1$ dimensional sphere in $\mathbb{R}^{m_1}$), $\omega_2 \in S^{m_2-1}$ and $u \in [\frac{\cosh s}{\cosh t}, 1]$. We have

$$\psi(t,s) = C \sinh t (\cosh t)^{m_2} \int_{\frac{\cosh s}{\cosh t}}^{1} (u \cosh t - \cosh s)^{(m_1-2)/2} (1-u^2)^{(m_2-2)/2} \, du,$$

which easily proves (2.5). The computation of the function $\psi$ is slightly easier if $m_2 = 0$ and the result is also given by (2.5). □

Our next proposition explains the role of the Lorentz space $L^{2,1}(G//K)$ which, by definition, is the subspace of $K$-bi-invariant functions in $L^{2,1}(G)$:

PROPOSITION 2. *The Abel transform*

$$\mathcal{A}f(a) = e^{\rho(\log a)} \int_N f(\overline{n}a) \, d\overline{n}$$

is bounded from $L^{2,1}(G//K)$ to $L^\infty(A/W)$. In other words, if $f$ is a locally integrable $K$-bi-invariant function on $G$ and $a \in A$ then:

$$e^{\rho(\log a)} \int_{\overline{N}} f(\overline{n}a)\, d\overline{n} \leq C ||f||_{L^{2,1}(G)}. \tag{2.7}$$

*Proof of Proposition* 2. The usual theory of Lorentz spaces (see, for example, [11, Chapter V]) shows that it suffices to prove the inequality (2.7) under the additional assumption that $f$ is the characteristic function of an open $K$-bi-invariant set of finite measure. For any $t \geq 0$, let $F(t) = f(a(t))$, so

$$||f||_{L^{2,1}(G)} = C \left[ \int_{\mathbb{R}_+} F(t)(\sinh t)^{m_1}(\sinh 2t)^{m_2}\, dt \right]^{1/2}. \tag{2.8}$$

In view of Lemma 1 and (2.8), it suffices to prove that for any $s \in \mathbb{R}$

$$\int_{t \geq |s|} F(t)\psi(t,s)\, dt \leq C \left[ \int_{\mathbb{R}_+} F(t)(\sinh t)^{m_1}(\sinh 2t)^{m_2}\, dt \right]^{1/2} \tag{2.9}$$

for any measurable function $F : \mathbb{R}_+ \to \{0,1\}$. Notice that if $t \geq 1 + |s|$ then $\psi(t,s) \approx e^{|\rho|t}$, $(\sinh t)^{m_1}(\sinh 2t)^{m_2} \approx e^{2|\rho|t}$ and it follows from Lemma 3 below that

$$\int_{t \geq |s|+1} F(t)\psi(t,s)\, dt \leq C \left[ \int_{t \geq |s|+1} F(t)(\sinh t)^{m_1}(\sinh 2t)^{m_2}\, dt \right]^{1/2}. \tag{2.10}$$

In order to deal with the integral in $t$ over the interval $[|s|, |s|+1]$ we consider two cases: $|s| \geq 1$ and $|s| \leq 1$. If $|s| \geq 1$ and $t \in [|s|, |s|+1]$, then $\psi(t,s) \approx e^{|\rho||s|}(t-|s|)^{(m_1+m_2-2)/2}$, $(\sinh t)^{m_1}(\sinh 2t)^{m_2} \approx e^{2|\rho||s|}$ and, since $(m_1+m_2-2)/2 \geq -1/2$, it follows that

$$\int_{|s|}^{|s|+1} F(t)\psi(t,s)\, dt \leq C e^{|\rho||s|} \int_{|s|}^{|s|+1} F(t)(t-|s|)^{-1/2}\, dt$$

$$= C e^{|\rho||s|} \int_0^1 F(|s|+u^2)\, du \leq C \left[ e^{2|\rho||s|} \int_0^1 F(|s|+u^2)u\, du \right]^{1/2}$$

$$\leq C \left[ \int_{|s|}^{|s|+1} F(t)(\sinh t)^{m_1}(\sinh 2t)^{m_2}\, dt \right]^{1/2}.$$

One of the inequalities in the sequence above follows from the estimate (2.11) below. This, together with (2.10), completes the proof of the proposition in the case $|s| \geq 1$. The estimation of the integrals over the interval $[|s|, |s|+1]$ is similar in the case $|s| \leq 1$. $\square$

LEMMA 3. *If $\delta \neq 0$ and $d\mu_1(t) = e^{\delta t}dt$, $d\mu_2(t) = e^{2\delta t}dt$ are two measures on $\mathbb{R}$ then*

$$||f||_{L^1(\mathbb{R}, d\mu_1)} \leq C_\delta ||f||_{L^{2,1}(\mathbb{R}, d\mu_2)}.$$

*Proof of Lemma* 3. One can assume that $f$ is the characteristic function of a set. The change of variable $t = (\log s)/\delta$ and the substitution $g(s) = f((\log s)/\delta)$ show that it suffices to prove that

$$\text{(2.11)} \qquad \frac{1}{|\delta|}\int_{\mathbb{R}_+} g(s)\,ds \leq C_\delta \left[\frac{1}{|\delta|}\int_{\mathbb{R}_+} g(s)s\,ds\right]^{1/2}$$

for any measurable function $g : \mathbb{R}_+ \to \{0,1\}$, which follows by a rearrangement argument. □

## 3. Proof of the maximal theorem

For any locally integrable function $f : X \to \mathbb{C}$ let

$$\text{(3.1)} \qquad \widetilde{\mathcal{M}_2}f(z) = \sup_{r \geq 1} \frac{1}{|B(z,r)|^{1/2}} \int_{B(z,r)} |f(z')|\,dz'.$$

Most of this section will be devoted to the proof of the following theorem:

THEOREM 4. *The operator $\widetilde{\mathcal{M}_2}$ is bounded from $L^{2,1}(X)$ to $L^{2,\infty}(X)$.*

Notice that Theorem B is an easy consequence of Theorem 4: let

$$\mathcal{M}_2^0 f(z) = \sup_{z \in B, r(B) \leq 1} \frac{1}{|B|}\int_B f(z')\,dz',$$

$$\mathcal{M}_2^1 f(z) = \sup_{z \in B, r(B) \geq 1} \frac{1}{|B|}\int_B f(z')\,dz',$$

where $r(B)$ is the radius of the ball $B$. We can assume that the Killing form on the Lie algebra $\mathfrak{g}$ is normalized such that $|H_0| = 1$. Let $o = \{K\}$ be the origin of the symmetric space $X$. Then the ball $B(o,r)$ is equal to the set of points $\{ka(t) \cdot o : k \in K, t \in [0,r)\}$ and one clearly has $|B(o,r)| \approx r^{m_1+m_2+1}$ if $r \leq 1$ and $|B(o,r)| \approx e^{2|\rho|r}$ if $r \geq 1$. The operator $\mathcal{M}_2^0$, the local part of $\mathcal{M}_2$, is clearly bounded on $L^p(X)$ for any $p > 1$. On the other hand, if $z$ belongs to a ball $B$ of radius $r \geq 1$, then $B(z,2r)$ contains the ball $B$ and $|B(z,2r)| \approx e^{2|\rho|\cdot 2r} \approx |B|^2$. Therefore

$$\frac{1}{|B|}\int_B f(z')\,dz' \leq \frac{C}{|B(z,2r)|^{1/2}}\int_{B(z,2r)} f(z')\,dz'$$

which shows that $\mathcal{M}_2^1 f(z) \leq C\widetilde{\mathcal{M}_2}f(z)$, and the conclusion of Theorem B follows by interpolation with the trivial $L^\infty$ estimate.

*Proof of Theorem* 4. Let $\chi_r$ be the characteristic function of the $K$-bi-invariant set $\{g \in G : d(g \cdot o, o) < r\}$. Since the measure of a ball of radius $r$ in $X$ is proportional to $e^{2|\rho|r}$ if $r \geq 1$, one has

$$\widetilde{\mathcal{M}_2}f(g \cdot o) \approx \sup_{r \geq 1}\left[e^{-|\rho|r}\int_G f(g' \cdot o)\chi_r(g'^{-1}g)\,dg'\right].$$

The change of variables $g = \overline{n}a(t)k$, $g' = \overline{n}'a(t')k'$ and the integral formula (2.4) show that

(3.2)
$$\widetilde{\mathcal{M}_2}f(\overline{n}a(t) \cdot o)$$
$$\leq C \sup_{r \geq 1}\left[e^{-|\rho|r}\int_{\mathbb{R}}\left(\int_{\overline{N}} f(\overline{n}'a(t') \cdot o)\chi_r(a(-t')\overline{n}'^{-1}\overline{n}a(t))\,d\overline{n}'\right)e^{2|\rho|t'}\,dt'\right].$$

We first deal with the integral over the group $\overline{N}$ and dominate the right-hand side of (3.2) using a standard maximal operator on the nilpotent group $\overline{N}$. For any $u > 0$ let $\mathcal{B}_u$ be the ball in $\overline{N}$ defined as the set $\{\overline{n}(v,w) : |v| \leq u$ and $|w| \leq u^2\}$. Clearly, $\int_{\mathcal{B}_u} 1\,d\overline{n} = Cu^{2|\rho|}$. The group $\overline{N}$ is equipped with non-isotropic dilations $\delta_u(\overline{n}(v,w)) = \overline{n}(uv, u^2 w)$, which are group automorphisms, therefore the maximal operator

$$\mathcal{N}g(\overline{n}) = \sup_{u>0}\left[\frac{1}{u^{2|\rho|}}\int_{\mathcal{B}_u}|g(\overline{n}\,\overline{m}^{-1})|\,d\overline{m}\right]$$

is bounded from $L^p(\overline{N})$ to $L^p(\overline{N})$ for any $p > 1$ ([13, Lemma 2.2]). For any locally integrable function $f : X \to \mathbb{R}_+$ and any $\overline{n} \in \overline{N}$ and $a \in A$ let

$$\mathcal{M}_3 f(\overline{n}a \cdot o) = \sup_{u>0}\left[\frac{1}{u^{2|\rho|}}\int_{\mathcal{B}_u}|f(\overline{n}\,\overline{m}^{-1}a \cdot o)|\,d\overline{m}\right].$$

Since the maximal operator $\mathcal{N}$ is bounded on $L^p(\overline{N})$ one has $\|\mathcal{M}_3 f\|_{L^p(X)} \leq C_p\|f\|_{L^p(X)}$ for any $p > 1$. We will now use the function $\mathcal{M}_3 f$ to control the integral over $\overline{N}$ in (3.2). Notice that (2.1) and (2.2), together with the fact that $d(ka(t) \cdot o, o) = t$ for any $t \geq 0$ and $k \in K$, show that if $\chi_r(a(-t')\overline{m}a(t)) = 1$ for some $\overline{m} \in \overline{N}$ then $\overline{m}$ has to belong to the ball $\mathcal{B}_{e^{(r-t-t')/2}}$; therefore

$$\int_{\overline{N}} f(\overline{n}'a(t') \cdot o)\chi_r(a(-t')\overline{n}'^{-1}\overline{n}a(t))\,d\overline{n}' \leq \int_{\mathcal{B}_{e^{(r-t-t')/2}}} f(\overline{n}\,\overline{m}^{-1}a(t') \cdot o)\,d\overline{m}$$
$$\leq Ce^{|\rho|(r-t-t')}\mathcal{M}_3 f(\overline{n}a(t') \cdot o).$$

If we substitute this inequality into (3.2) we conclude that

(3.3) $$\widetilde{\mathcal{M}_2}f(\overline{n}a(t) \cdot o) \leq Ce^{-|\rho|t}\int_{\mathbb{R}} \mathcal{M}_3 f(\overline{n}a(t') \cdot o)e^{|\rho|t'}\,dt'.$$

We can now estimate the $L^{2,\infty}$ norm of $\widetilde{\mathcal{M}_2}f$: for any $\lambda > 0$, the set $E_\lambda = \{z \in X : \widetilde{\mathcal{M}_2}f(z) > \lambda\}$ is included in the set

$$\{\overline{n}a(t) \cdot o : e^{-|\rho|t}\int_{\mathbb{R}} \mathcal{M}_3 f(\overline{n}a(t') \cdot o)e^{|\rho|t'}\,dt' > \lambda/C\}.$$

The measure $dz$ in $X$ is proportional to the measure $e^{2|\rho|t}\,d\overline{n}\,dt$ in $\overline{N}\times\mathbb{R}$ under the identification $z = \overline{n}a(t)\cdot o$. Therefore the measure of this last set is less than or equal to

$$\frac{C\int_{\overline{N}}\left[\int_{\mathbb{R}}\mathcal{M}_3 f(\overline{n}a(t')\cdot o)e^{|\rho|t'}\,dt'\right]^2 d\overline{n}}{\lambda^2};$$

hence

$$(3.4) \qquad \|\widetilde{\mathcal{M}_2}f\|_{L^{2,\infty}}^2 \leq C\int_{\overline{N}}\left[\int_{\mathbb{R}}\mathcal{M}_3 f(\overline{n}a(t')\cdot o)e^{|\rho|t'}\,dt'\right]^2 d\overline{n}.$$

One can now use the following simple lemma to dominate the right-hand side of (3.4):

LEMMA 5. *If $U$ and $V$ are two measure spaces with measures $du$ and $dv$ respectively, and $H: U\times V \to \mathbb{R}_+$ is measurable then*

$$\left[\int_U \|H(u,.)\|_{L^{2,1}(V,dv)}^2\,du\right]^{1/2} \leq C\|H\|_{L^{2,1}(U\times V,\,du\,dv)}.$$

The proof of this lemma is straightforward. Combining Lemma 3 (at the end of the previous section) and Lemma 5, one has

(3.5)
$$\begin{aligned}
\int_{\overline{N}}\left[\int_{\mathbb{R}}\mathcal{M}_3 f(\overline{n}a(t')\cdot o)e^{|\rho|t'}\,dt'\right]^2 d\overline{n} &\leq C\int_{\overline{N}}\|\mathcal{M}_3 f(\overline{n}a(.)\cdot o)\|_{L^{2,1}(\mathbb{R},e^{2|\rho|t'}\,dt')}^2\,d\overline{n} \\
&\leq C\|\mathcal{M}_3 f(\overline{n}a(t')\cdot o)\|_{L^{2,1}(\overline{N}\times\mathbb{R},e^{2|\rho|t'}\,d\overline{n}\,dt')}^2 \\
&\leq C\|\mathcal{M}_3 f\|_{L^{2,1}(X)}^2.
\end{aligned}$$

Finally, since the maximal operator $\mathcal{M}_3$ is bounded on $L^p(X)$ for any $p>1$, it follows by the general version of Marcinkiewicz interpolation theorem that $\|\mathcal{M}_3 f\|_{L^{2,1}(X)} \leq C\|f\|_{L^{2,1}(X)}$ and Theorem 4 follows from (3.4) and (3.5). □

## 4. A covering lemma

A simple connection between covering lemmas and boundedness of maximal operators is explained in [2]. In our setting we have:

COROLLARY 6. *If a collection of balls $B_i \subset X$, $i\in I$, has the property that $|\cup B_i| < \infty$ then one can select a finite subset $J \subset I$ such that*

$$(4.1) \qquad \text{(i)} \qquad \left|\bigcup_{i\in I} B_i\right| \leq C\left|\bigcup_{j\in J} B_j\right|;$$

$$\text{(ii)} \qquad \left\|\sum_{j\in J}\chi_{B_j}\right\|_{L^{2,\infty}(X)} \leq C\left|\bigcup_{i\in I} B_i\right|^{1/2}.$$

It follows from (4.1) that

$$\left\|\sum_{j\in J}\chi_{B_j}\right\|_{L^q(X)} \leq C_q \left|\bigcup_{i\in I} B_i\right|^{1/q}$$

for any $q \in [1, 2)$. Thus, in the terminology of [2], the family of natural balls on symmetric spaces of real rank one has the covering property $V_q$ if and only if $q \in [1, 2)$.

## 5. Proof of the convolution theorem

In this section we will prove Theorem A. In view of the general theory of Lorentz spaces, it suffices to prove that

$$(5.1) \qquad \iint_{G\times G} f(z)g(z^{-1}z')h(z')\,dz'\,dz \leq C\|f\|_{L^{2,1}}\|g\|_{L^{2,1}}\|h\|_{L^{2,1}}$$

whenever $f, g, h : G \to \{0, 1\}$ are characteristic functions of open sets of finite measure. We can also assume that $g$ is supported away from the origin of the group, for example in the set $\bigcup_{t>1} Ka(t)K$. The main part of our argument is devoted to proving that the left-hand side of (5.1) is controlled by an integral involving suitable rearrangements of the functions $f$, $g$ and $h$, as in (5.19). Let $z = \bar{n}a(t)k$, $z' = \bar{n}'a(t')k'$ and the left-hand side of (5.1) becomes

$$(5.2) \qquad \int_K \int_K \int_{\mathbb{R}} \int_{\mathbb{R}} I(k,k',t,t')e^{2|\rho|(t+t')}\,dt'\,dt\,dk'\,dk,$$

where

(5.3)
$$I(k,k',t,t') = \iint_{\overline{N}\times\overline{N}} f(\bar{n}a(t)k)g(k^{-1}a(-t)\bar{n}^{-1}\bar{n}'a(t')k')h(\bar{n}'a(t')k')\,d\bar{n}'\,d\bar{n}$$

We will show how to dominate the expression in (5.2) in four steps.

*Step* 1. *Integration on the subgroup* $\overline{N}$. As in the proof of the maximal theorems, we start by integrating on $\overline{N}$. Define $F_1, H_1 : K \times \mathbb{R} \to \mathbb{R}_+$ by

$$F_1(k,t) = \int_{\overline{N}} f(\bar{n}a(t)k)\,d\bar{n}$$

and

$$H_1(k',t') = \int_{\overline{N}} h(\bar{n}'a(t')k')\,d\bar{n}'.$$

Using the simple inequality

$$\iint_{\overline{N}\times\overline{N}} a(\bar{n})b(\bar{n}^{-1}\bar{n}')c(\bar{n}')\,d\bar{n}'\,d\bar{n}$$
$$\leq \left(\int_{\overline{N}} b(\bar{n})\,d\bar{n}\right)\left[\min\left(\left(\int_{\overline{N}} a(\bar{n})\,d\bar{n}\right), \left(\int_{\overline{N}} c(\bar{n})\,d\bar{n}\right)\right)\right],$$

which holds for any measurable functions $a, b, c : \overline{N} \to [0,1]$ with compact support, it follows that the integral $I(k, k', t, t')$ in (5.3) is dominated by

$$(5.4) \qquad \min\left[F_1(k,t), H_1(k', t')\right] \left[\int_{\overline{N}} g(k^{-1} a(-t) \overline{n}_1 a(t') k') \, d\overline{n}_1\right].$$

By (2.2), the map $\overline{n}_1 \to a(-t)\overline{n}_1 a(t) = \overline{n}_2$ is a dilation of $\overline{N}$ with $d\overline{n}_1 = e^{-2|\rho|t} d\overline{n}_2$; therefore

$$(5.5) \quad \int_{\overline{N}} g(k^{-1} a(-t) \overline{n}_1 a(t') k') \, d\overline{n}_1 = e^{-2|\rho|t} \int_{\overline{N}} g(k^{-1} \overline{n}_2 a(t'-t) k') \, d\overline{n}_2$$

$$= C e^{-2|\rho|t} \int_{\mathbb{R}^{m_1} \times \mathbb{R}^{m_2}} g(k^{-1} \overline{n}(v, w) a(t'-t) k') \, dv \, dw$$

$$= C e^{-2|\rho|t} \int_{u \geq |t'-t|} \int_{T_{u, t'-t}} g(k^{-1} \overline{n}(P) a(t'-t) k') \, d\omega_{u, t'-t}(P) \, du.$$

The surfaces $T_{u,s}$ defined in (2.6) for $\{(u,s) \in \mathbb{R}_+ \times \mathbb{R} : u \geq |s|\}$ and the associated measures $d\omega_{u,s}$ have the same meaning as in the proof of Lemma 1. Let

$$(5.6) \quad G_1(k, k', u, s) = \left(\int_{T_{u,s}} 1 \, d\omega_{u,s}\right)^{-1} \left[\int_{T_{u,s}} g(k^{-1} \overline{n}(P) a(s) k') \, d\omega_{u,s}(P)\right]$$

be the average of the function $P \to g(k^{-1} \overline{n}(P) a(s) k')$ on the surface $T_{u,s}$ (the domain of definition of $G_1$ is $\{(k, k', u, s) \in K \times K \times \mathbb{R}_+ \times \mathbb{R} : u \geq |s|\}$, and $G_1(k, k', u, s) \in [0, 1]$). If we substitute this definition in (5.5), we conclude that

$$\int_{\overline{N}} g(k^{-1} a(-t) \overline{n}_1 a(t') k') \, d\overline{n}_1$$
$$= C e^{-|\rho|(t+t')} \int_{u \geq |t'-t|} G_1(k, k', u, t'-t) \psi(u, t'-t) \, du.$$

The function $\psi(u, s)$ was computed in the proof of Lemma 1 and is given by (2.5). Finally, if we substitute this last formula in (5.4), we find that the integral $I(k, k', t, t')$ is dominated by

$$C e^{-|\rho|(t+t')} \min\left[F_1(k,t), H_1(k', t')\right] \int_{u \geq |t'-t|} G_1(k, k', u, t'-t) \psi(u, t'-t) \, du,$$

which shows that the left-hand side of (5.1) is dominated by

$$(5.7) \qquad C \int_K \int_K \int_{\mathbb{R}} \int_{\mathbb{R}} \int_{u \geq |t'-t|} \min\left[F_1(k,t), H_1(k', t')\right]$$
$$G_1(k, k', u, t'-t) \psi(u, t'-t) e^{|\rho|(t+t')} \, du \, dt' \, dt \, dk' \, dk.$$

For later use, we record the following properties of the functions $F_1$ and $H_1$:

$$(5.8) \qquad \begin{aligned} \|f\|_{L^{2,1}(G)} &= \left[C_2 \int_K \int_{\mathbb{R}} F_1(k,t) e^{2|\rho|t} \, dt \, dk\right]^{1/2}, \\ \|h\|_{L^{2,1}(G)} &= \left[C_2 \int_K \int_{\mathbb{R}} H_1(k', t') e^{2|\rho|t'} \, dt' \, dk'\right]^{1/2}. \end{aligned}$$

*Step* 2. *Integration on the subgroup A.* Let $\chi_1$ and $\chi_2$, be the characteristic functions of the sets $\{(k, k', t, t') : F_1(k, t) \leq H_1(k', t')\}$ and $\{(k, k', t, t') : H_1(k', t') \leq F_1(k, t)\}$ respectively. For any $k, k', t, t'$ one has

(5.9)
$$\begin{cases} F_1(k, t)\chi_1(k, k', t, t') \leq H_1(k', t'), \\ H_1(k', t')\chi_2(k, k', t, t') \leq F_1(k, t). \end{cases}$$

Since $\chi_1 + \chi_2 \geq 1$, the expression (5.7) is less than or equal to the sum of two similar expressions of the form

$$C \int_K \int_K \int_{\mathbb{R}} \int_{\mathbb{R}} \int_{u \geq |t'-t|} F_1(k, t)\chi_1(k, k', t, t')$$
$$G_1(k, k', u, t' - t)\psi(u, t' - t)e^{|\rho|(t+t')} \, du \, dt' \, dt \, dk' \, dk.$$

The change of variable $t' = t + s$ in the expression above shows that it is equal to

(5.10)
$$C \int_K \int_K \int_{\mathbb{R}} \int_{\mathbb{R}} \int_{u \geq |s|} F_1(k, t)\chi_1(k, k', t, t+s)$$
$$G_1(k, k', u, s)\psi(u, s)e^{2|\rho|t}e^{|\rho|s} \, du \, dt \, ds \, dk' \, dk,$$

and the first of the inequalities in (5.9) becomes

(5.11)
$$F_1(k, t)\chi_1(k, k', t, t+s)e^{2|\rho|t} \leq H_1(k', t+s)e^{2|\rho|t}.$$

Let $F(k) = \left[\int_{\mathbb{R}} F_1(k, t)e^{2|\rho|t} \, dt\right]^{1/2}$, $H(k') = \left[\int_{\mathbb{R}} H_1(k', t')e^{2|\rho|t'} \, dt'\right]^{1/2}$ and

$$A(k, k', s) = \int_{\mathbb{R}} F_1(k, t)\chi_1(k, k', t, t+s)e^{2|\rho|t} \, dt.$$

The expression (5.10) becomes

(5.12)
$$C \int_K \int_K \int_{\mathbb{R}} \int_{u \geq |s|} A(k, k', s)G_1(k, k', u, s)\psi(u, s)e^{|\rho|s} \, du \, ds \, dk' \, dk.$$

Clearly, $A(k, k', s) \leq F(k)^2$ (since $\chi_1 \leq 1$) and $A(k, k', s) \leq e^{-2|\rho|s}H(k')^2$ by (5.11); therefore

$$e^{|\rho|s}A(k, k', s) \leq \begin{cases} e^{|\rho|s}F(k)^2 & \text{if } e^{|\rho|s} \leq H(k')/F(k), \\ e^{-|\rho|s}H(k')^2 & \text{if } e^{|\rho|s} \geq H(k')/F(k). \end{cases}$$

If we substitute this inequality in (5.12) we find that the left-hand side of (5.1) is dominated by

(5.13)
$$C \int_K \int_K \int_{e^{|\rho|s} \leq H(k')/F(k)} \int_{u \geq |s|} F(k)^2 G_1(k, k', u, s)\psi(u, s)e^{|\rho|s} \, du \, ds \, dk' \, dk$$
$$+ C \int_K \int_K \int_{e^{|\rho|s} \geq H(k')/F(k)} \int_{u \geq |s|} H(k')^2 G_1(k, k', u, s)\psi(u, s)e^{-|\rho|s} \, du \, ds \, dk' \, dk.$$

We pause for a moment to note that our estimates so far, together with the proof of Lemma 1 in the second section, suffice to prove that $L^{2,1}(G) * L^{2,1}(G//K) \subseteq L^{2,\infty}(G)$: if $g$ is a $K$-bi-invariant function, then $G_1(k, k', u, s)$ depends only on $u$, and (2.9) shows that

$$\int_{u \geq |s|} G_1(k, k', u, s) \psi(u, s) \, du \leq C \|g\|_{L^{2,1}}.$$

As a consequence, both terms in (5.13) are dominated by

$$C \|g\|_{L^{2,1}} \int_K \int_K F(k) H(k') \, dk' \, dk;$$

therefore

$$\iint_{G \times G} f(z) g(z^{-1} z') h(z') \, dz' \, dz \leq C \|g\|_{L^{2,1}} \int_K \int_K F(k) H(k') \, dk' \, dk$$
$$\leq C \|f\|_{L^{2,1}} \|g\|_{L^{2,1}} \|h\|_{L^{2,1}}.$$

Here we used the fact that, as a consequence of (5.8),

(5.14)
$$\|f\|_{L^{2,1}(G)} = \left[ C_2 \int_K F(k)^2 \, dk \right]^{1/2},$$
$$\|h\|_{L^{2,1}(G)} = \left[ C_2 \int_K H(k')^2 \, dk' \right]^{1/2}.$$

*Step* 3. *A rearrangement inequality.* In the general case (if $g$ is not assumed to be $K$-bi-invariant) we will show that both terms in (5.13) are dominated by some expression of the form

$$C \int_0^1 \int_0^1 \int_{\mathbb{R}_+} F^*(x) H^*(y) G^{**}(x, y, u) e^{|\rho|u} \, du \, dy \, dx$$

where $F^*, H^* : (0, 1] \to \mathbb{R}_+$ are the usual nonincreasing rearrangements of the functions $F$ and $H$ (recall that the measure of $K$ is equal to 1) and $G^{**} : (0, 1] \times (0, 1] \times \mathbb{R}_+ \to \{0, 1\}$ is a suitable "double" rearrangement of $g$. The precise definitions are the following: if $a : K \to \mathbb{R}_+$ is a measurable function then the nonincreasing rearrangement $a^* : (0, 1] \to \mathbb{R}_+$ is the right semicontinuous nonincreasing function with the property that

$$|\{k \in K : a(k) > \lambda\}| = |\{x \in (0, 1] : a^*(x) > \lambda\}| \text{ for any } \lambda \in [0, \infty).$$

Assume now that $a : K \times K \to \mathbb{R}_+$ is a measurable function. For almost every $k \in K$ let $a^*(k, y)$, $y \in (0, 1]$, be the nonincreasing rearrangement of the function $k' \to a(k, k')$ and let $a^{**}(x, y)$ be the nonincreasing rearrangement of the function $k \to a(k, y)$ (clearly $a^{**} : (0, 1] \times (0, 1] \to \mathbb{R}_+$). The following lemma summarizes some of the well-known properties of nonincreasing rearrangements (see for example [11, Chapter V]):

LEMMA 7. (a) *If $a : K \to \mathbb{R}_+$ is a measurable function then*

$$\left[\int_K a(k)^2 \, dk\right]^{1/2} = \left[\int_{(0,1]} a^*(x)^2 \, dx\right]^{1/2}.$$

(b) *If $a : K \times K \to \mathbb{R}_+$ is a measurable function then*

(i) $$\int_K \int_K a(k, k') \, dk' \, dk = \int_0^1 \int_0^1 a^{**}(x, y) \, dy \, dx.$$

(ii) *The function $a^{**}$ is nonincreasing: $a^{**}(x, y) \leq a^{**}(x', y')$ whenever $x \geq x'$ and $y \geq y'$.*

(iii) *For any measurable sets $D, E \subset K$ with measures $|D|$ and $|E|$*

$$\int_D \int_E a(k, k') \, dk' \, dk \leq \int_0^{|D|} \int_0^{|E|} a^{**}(x, y) \, dy \, dx.$$

Returning to our setting, let $F^*$ and $H^*$ be the nonincreasing rearrangements of $F$ and $H$, let $\tilde{g} : K \times K \times \mathbb{R}_+ \to \{0, 1\}$ be given by $\tilde{g}(k, k', u) = g(k^{-1}a(u)k')$ and let $G^{**} : (0, 1] \times (0, 1] \times \mathbb{R}_+ \to \{0, 1\}$ be the double rearrangement of the function $\tilde{g}$ (i.e., $G^{**}(.,.,u)$ is the double rearrangement of $\tilde{g}(.,.,u)$ for all $u \geq 0$). Recall that we assumed that the function $g$ is the characteristic function of a set included in $\bigcup_{u > 1} Ka(u)K$; therefore

$$(5.15) \qquad ||g||_{L^{2,1}(G)} \approx \left[\int_{\mathbb{R}_+} \int_0^1 \int_0^1 G^{**}(x, y, u) e^{2|\rho|u} \, dy \, dx \, du\right]^{1/2}.$$

We will now show how to use these rearrangements to dominate the two expressions in (5.13). For any integers $m$, $n$ let $D_m = \{k \in K : F(k) \in [e^{|\rho|m}, e^{|\rho|(m+1)}]\}$, $E_n = \{k' \in K : H(k') \in [e^{|\rho|n}, e^{|\rho|(n+1)}]\}$ and let $D_{-\infty} = \{k \in K : F(k) = 0\}$, $E_{-\infty} = \{k' \in K : H(k') = 0\}$ such that $K = \bigcup_m D_m = \bigcup_n E_n$. Let $\delta_m$, respectively $\varepsilon_n$, be the measures of the sets $D_m$, respectively $E_n$, as subsets of $K$. The first of the two expressions in (5.13) is dominated by

(5.16)
$$C \sum_{m,n} \int_{D_m} \int_{E_n} \int_{s \leq (n-m+1)} \int_{u \geq |s|} e^{2|\rho|(m+1)} G_1(k, k', u, s) \psi(u, s) e^{|\rho|s} \, du \, ds \, dk' \, dk.$$

Combining the definition (5.6) of the function $G_1$ (recall that the surfaces $T_{u,s}$ are defined as the set of points $P \in \mathbb{R}^{m_1} \times \mathbb{R}^{m_2}$ with the property that $\overline{n}(P)a(s) \in Ka(u)K$), the fact that $dk$ is a Haar measure on $K$ and the last statement of Lemma 7, we conclude that

$$\int_{D_m} \int_{E_n} G_1(k, k', u, s) \, dk' \, dk \leq \int_0^{\delta_m} \int_0^{\varepsilon_n} G^{**}(x, y, u) \, dy \, dx$$

for any $s$ with the property that $|s| \leq u$. Substituting this inequality in (5.16), we find that the expression in (5.16) is dominated by

$$(5.17) \quad C \sum_{m,n} \int_{\mathbb{R}_+} e^{2|\rho|m} \left[ \int_0^{\delta_m} \int_0^{\varepsilon_n} G^{**}(x,y,u)\,dy\,dx \right]$$

$$\left[ \int_{s \leq (n-m+1), |s| \leq u} \psi(u,s) e^{|\rho|s}\,ds \right] du.$$

The formula (2.5) shows that the last of the integrals in the expression above is dominated by $Ce^{|\rho|u}e^{|\rho|(n-m)}$; therefore the first of the two expressions in (5.13) is dominated by

$$(5.18) \quad C \int_{\mathbb{R}_+} \sum_{m,n} \left[ e^{|\rho|(m+n)} \int_0^{\delta_m} \int_0^{\varepsilon_n} G^{**}(x,y,u)\,dy\,dx \right] e^{|\rho|u}\,du.$$

Let

$$S(x,y) = \sum_{m,n} \left[ e^{|\rho|(m+n)} \chi_{\delta_m}(x) \chi_{\varepsilon_n}(y) \right],$$

where $\chi_{\delta_m}$, $\chi_{\varepsilon_n}$ are the characteristic functions of sets $(0, \delta_m)$, respectively $(0, \varepsilon_n)$. If $m_x = \max\{m : \delta_m > x\}$ and $n_y = \max\{n : \varepsilon_n > y\}$ then $S(x,y) \leq Ce^{|\rho|(m_x+n_y)}$. Clearly $F^*(x) \geq e^{|\rho|m_x}$, $H^*(y) \geq e^{|\rho|n_y}$; therefore the expression (5.18) is dominated by

$$C \int_{\mathbb{R}_+} \int_0^1 \int_0^1 F^*(x) H^*(y) G^{**}(x,y,u) e^{|\rho|u}\,dy\,dx\,du.$$

One can deal with the second of the two expressions in (5.13) in a similar way; therefore

$$(5.19) \quad \iint_{G \times G} f(z) g(z^{-1}z') h(z')\,dz'\,dz$$

$$\leq C \int_{\mathbb{R}_+} \int_0^1 \int_0^1 F^*(x) H^*(y) G^{**}(x,y,u) e^{|\rho|u}\,dy\,dx\,du.$$

*Step* 4. *Final estimates.* Let $\mathcal{K}$ be a suitable constant (to be chosen later) and let $\mathcal{U} = \{(x,y,u) : F^*(x) H^*(y) \leq \mathcal{K} e^{|\rho|u}\}$ and $\mathcal{V} = \{(x,y,u) : F^*(x) H^*(y) \geq \mathcal{K} e^{|\rho|u}\}$. By (5.15),

$$\int_{\mathcal{U}} F^*(x) H^*(y) G^{**}(x,y,u) e^{|\rho|u}\,dy\,dx\,du$$

$$\leq \int_{\mathbb{R}_+} \int_0^1 \int_0^1 \mathcal{K} G^{**}(x,y,u) e^{2|\rho|u}\,dy\,dx\,du \leq C\mathcal{K} \|g\|_{L^{2,1}}^2.$$

Using Lemma 7(a), (5.14) and the fact that $G^{**}(x,y,u) \leq 1$ one has

$$\int_{\mathcal{V}} F^*(x)H^*(y)G^{**}(x,y,u)e^{|\rho|u}\,dy\,dx\,du \;\leq\; C\int_0^1\int_0^1 \frac{[F^*(x)H^*(y)]^2}{\mathcal{K}}\,dy\,dx$$

$$\leq\; C\frac{||f||^2_{L^{2,1}}||h||^2_{L^{2,1}}}{\mathcal{K}}.$$

Finally one lets $\mathcal{K} = (||g||_{L^{2,1}})^{-1}(||f||_{L^{2,1}}||h||_{L^{2,1}})$ and the theorem follows.

## 6. A general rearrangement inequality

We will now extend the rearrangement inequality (5.19) to the case when $f$, $g$, $h$ are arbitrary measurable functions (not just characteristic functions of sets). For any measurable function $f : G \to \mathbb{R}_+$ we define the function $F^* : (0,1] \to \mathbb{R}_+$ by the following procedure: first, let $\tilde{f} : K \times (0,\infty) \to \mathbb{R}_+$ be defined, for almost every $k \in K$, as the usual nonincreasing rearrangement of the function $f_k : \overline{N} \times A \to \mathbb{R}_+$, $f_k(\overline{n}a) = f(\overline{n}ak)$ with respect to the measure $e^{2\rho(\log a)}d\overline{n}da$. Using the function $\tilde{f}$ we define the function $\tilde{F} : (0,1] \times (0,\infty) \to \mathbb{R}_+$: for each $r > 0$ fixed, the function $\tilde{F}(.,r)$ is the usual the nonincreasing rearrangement of the function $k \to \tilde{f}(k,r)$. Finally let

(6.1) $$F^*(x) = \frac{1}{2}\int_0^\infty \tilde{F}(x,r)r^{-1/2}\,dr$$

be the $L^{2,1}$ norm of the function $r \to \tilde{F}(x,r)$. Notice that this definition of the function $F^*$ agrees with our earlier definition if $f$ is a characteristic function.

THEOREM 8. *If $f, g, h : G \to R_+$ are measurable functions then*

(6.2) $$\iint_{G\times G} f(z)g(z^{-1}z')h(z')\,dz'\,dz$$
$$\leq\; C\int_{\mathbb{R}_+}\int_0^1\int_0^1 F^*(x)H^*(y)G^{**}(x,y,u)\phi(u)\,dy\,dx\,du,$$

*where $G^{**} : (0,1] \times (0,1] \times \mathbb{R}_+ \to \mathbb{R}_+$ is the double rearrangement of the function $(k,k',u) \to g(k^{-1}a(u)k')$ (the same definition as before), $F^*$ and $H^*$ are defined in the previous paragraph and $\phi(u) = u^{m_1+m_2}$ if $u \leq 1$ and $\phi(u) = e^{|\rho|u}$ if $u \geq 1$.*

*Proof of Theorem 8.* Notice that

$$\phi(u) \approx \sup_{r\in[-u,u]} e^{-|\rho|r}\int_{s\leq r, |s|\leq u} \psi(u,s)e^{|\rho|s}\,ds.$$

Notice also that if $f$ and $h$ are characteristic functions of sets then (6.2) is equivalent to (5.19). If $f, h$ are simple positive functions, one can write (uniquely up to sets of measure zero) $f = \sum_1^{M_1} c_i f_i$, $h = \sum_1^{M_2} d_j h_j$, where $c_i, d_j > 0$ and $f_i$ and $h_j$, are characteristic functions of sets $U_i$ and $V_j$ with the property that for all $i$ and $j$ one has $U_{i+1} \subset U_i$ and $V_{j+1} \subset V_j$. Simple manipulations involving rearrangements show that $F^* = \sum_1^{M_1} c_i F_i^*$ and $H^* = \sum_1^{M_2} d_j H_j^*$ (this explains the reason why we chose the apparently complicated definition of the function $F^*$ in (6.1)), and (6.2) follows by summation. Finally, a standard argument shows that (6.2) holds for arbitrary measurable functions $f$, $g$ and $h$ for which the right-hand side integral in (6.2) converges. □

PRINCETON UNIVERSITY, PRINCETON, NJ
*Current address*: INSTITUTE FOR ADVANCED STUDY, PRINCETON, NJ
*E-mail address*: aionescu@math.ias.edu